\documentclass{gtmon_a}
\pdfoutput=1

\usepackage{amscd}
\usepackage[v2,cmtip]{xy}
\usepackage[T1,T5]{fontenc}

\proceedingstitle{Proceedings of the School and Conference in Algebraic
Topology (The Vietnam National University, Hanoi, 9-20 August 2004)}
\conferencestart{9 August 2004}
\conferenceend{20 August 2004}
\conferencename{School and Conference in Algebraic Topology}
\conferencelocation{Vietnam National University, Hanoi, Vietnam}

\editor{John Hubbuck}
\givenname{John}
\surname{Hubbuck}

\editor{Nguy\~\ecircumflex{}n H V H\uhorn{}ng}
\givenname{H\uhorn{}ng}
\surname{Nguy\~\ecircumflex{}n}

\editor{Lionel Schwartz}
\givenname{Lionel}
\surname{Schwartz}

\title{Sub-Hopf algebras of the Steenrod algebra\\and the Singer transfer}

\author{L\^e Minh H\`a}
\surname{L\^e}
\givenname{Minh H\`a}
\address{Department of Mathematics\\
Vietnam National University\\\newline
334 Nguyen Trai Street\\
Hanoi\\
Vietnam}
\email{minhha@vnu.edu.vn}
\email{leminhha.vnu@gmail.com}
\urladdr{}

\dedicatory{Dedicated to Professor Hu\`ynh M\`ui on the occasion of his
60th birthday}

\volumenumber{11}
\issuenumber{}
\publicationyear{2007}
\papernumber{5}
\startpage{81}
\endpage{105}

\doi{}
\MR{}
\Zbl{}

\arxivreference{}


\received{1 April 2005}
\revised{17 December 2005}
\accepted{2 March 2006}
\published{}
\publishedonline{14 November 2007}
\proposed{14 November 2007}
\seconded{}
\corresponding{}
\editor{}
\version{}

\let\xysavmatrix\xymatrix
\def\xymatrix{\disablesubscriptcorrection\xysavmatrix}
\AtBeginDocument{\let\bar\wbar\let\tilde\wtilde\let\hat\what}

\def\co{\colon \thinspace}
\def\rmap#1{\, \smash{\mathop{\longrightarrow}\limits^{#1}} \:}

\def\f2{{\mathbb F}_{2}}
\def\z2{{\mathbb Z}/2}

\DeclareMathOperator{\st}{\mathcal{A}}
\DeclareMathOperator{\sst}{\mathcal{A}^{*}}

\DeclareMathOperator{\Hom}{Hom}

\def\invlim#1{\mathchoice
        {\mathop{\displaystyle\mathop{\mathrm{\rm lim}}_{\longleftarrow}}
        _{#1}}
        {\mathop{\smash{\displaystyle\mathop{\mathrm{\rm lim}}
        _{\longleftarrow}}}
        _{#1}\vphantom{\displaystyle\mathop{\mathrm{\rm lim}}
        _{\longleftarrow}}}
        {}{}}

\makeatletter
\def\cnewtheorem#1[#2]#3{\newtheorem{#1}{#3}[section]
\expandafter\let\csname c@#1\endcsname\c@theorem}

\newtheorem{theorem}{Theorem}[section]
\cnewtheorem{lemma}[theorem]{Lemma}
\cnewtheorem{proposition}[theorem]{Proposition}
\cnewtheorem{corollary}[theorem]{Corollary}
\cnewtheorem{conjecture}[theorem]{Conjecture}
\theoremstyle{definition}
\cnewtheorem{definition}[theorem]{Definition}
\cnewtheorem{example}[theorem]{Example}
\cnewtheorem{remark}[theorem]{Remark}
\cnewtheorem{question}[theorem]{Question}
\cnewtheorem{notation}[theorem]{Notation}

\makeatother  

\newcommand{\otigls}{\rlap{$\,\,\otimes$}\lower 7pt\hbox{$_{\gls}$}}
\newcommand{\otigl}{\rlap{$\,\,\otimes$}\lower 7pt\hbox{$\,_{GL}$}}
\newcommand{\otigln}[1]{\rlap{$\,\,\otimes$}\lower 7pt\hbox{$\,_{GL_{#1}}$}}
\newcommand{\oticala}{\,\rlap{$\otimes$}\lower 7pt\hbox{$_{\mathcal{A}}$}\,}
\newcommand{\doplus}[1]{\rlap{$\,\oplus\,$}\lower 7pt \hbox{$_{#1}$}}

\numberwithin{equation}{section}

\makeop{res}
\newcommand{\dslash}{{/}\!{/}}

\begin{document}

\begin{abstract}
The Singer transfer provides an interesting connection between modular
representation theory and the cohomology of the Steenrod algebra. We
discuss a version of ``Quillen stratification'' theorem for the Singer
transfer and its consequences.
\end{abstract}

\maketitle

\section{Introduction}\label{sec: intro}
Let $\st$ denote the mod $2$ Steenrod algebra
(see Steenrod and Epstein \cite{Steenrod-book}). The problem of computing
its cohomology
$H^{*,*} (\st)$  is of great importance in algebraic topology, for
this bigraded commutative algebra is the $E^{2}$ term  of the
Adams spectral sequence (see Adams \cite{Adams-CMH58}) converging to the
stable homotopy groups of spheres. But despite intensive investigation for nearly
half a century, the structure of this cohomology algebra remains
elusive. In fact, only recently was a complete description of generators and relations in
cohomological dimension $4$ given, by Lin and Mahowald 
\cite{Lin-Mahowald,Lin-TAMS}. In higher degrees, several infinite dimensional subalgebras of $H^{*}
(\st )$ have been constructed and studied. The first such subalgebra
\cite{Adams-CMH58}, called the Adams subalgebra, is generated by the
elements $h_{i} \in H^{1, 2^{i}} (\st )$ for $i \geq 0$. 
Mahowald and Tangora \cite{Mahowald-Tangora} constructed the so-called
\emph{wedge}
subalgebra which consists of some basic generators,
propagated by the Adams periodicity operator $P^{1}$ and by
multiplication with a certain element $g \in H^{4,20} (\st
)$\footnote{This element is also denoted as $g_{1}$ in literature.}.  The wedge subalgebra was subsequently
expanded by another kind of periodicity operator $M$, discovered by
Margolis, Priddy and Tangora \cite{Margolis-Priddy-Tangora}.  On the other hand, perhaps the most important
result to date on the structure of $H^{*} (\st )$ is a beautiful theorem of Palmieri~\cite{Palmieri99}
which gives a version of the  famous Quillen stratification theorem in
group cohomology for the  cohomology of the Steenrod algebra. Loosely speaking, Palmieri's
theorem says that modulo nilpotent elements, the
cohomology of the Steenrod algebra is completely determined by the
cohomology of its elementary sub-Hopf algebras. The underlying ideas  in both type of
results mentioned above are the same: One can obtain information
about $H^{*,*} (\st )$ by studying its restriction to various suitably
chosen sub-Hopf algebras of $\mathcal{A}$.

In this paper, we take up this idea to investigate the Singer transfer
\cite{Singer-transfer}.  To describe this map, we first need some
notation. We work over the field $\f2$ and let $V_{n}$ denote the
elementary abelian $2$--group of rank $n$. It is well-known that the
mod $2$ homology $H_{*} (BV_{n})$ is a divided power algebra on $n$
generators.  Furthermore, there is an action of the group algebra $\st
[GL(n, \f2)]$, where the Steenrod algebra acts by dualizing the canonical
action in cohomology, and the general linear group $GL(n):=GL(n,
\f2)$ acts by matrix substitution. Let $\wwbar{\mathcal{A}}$ be
the augmentation ideal of $\st$. Let $P_{\st} H_{*} (BV_{n})$ be
the subring of $H_{*} (BV_{n})$ consisting of all elements that are
$\wwbar{\mathcal{A}}$--annihilated. In \cite{Singer-transfer}, Singer
constructed a  map from this subring to the cohomology of
the Steenrod algebra
\[
\mathtt{Tr}_{n}^{\st }\co P_{\mathcal{A}} H_{*} (BV_{n}) \rmap{} H^{n,n+*} (\mathcal{A}),
\]
in such a way that the \textit{total} transfer $\mathtt{Tr}^{\st } =
\bigoplus_{n} \mathtt{Tr}_{n}^{\st }$ is a bigraded algebra homomorphism
with respect to the product by concatenation in the domain and the usual
Yoneda product for the Ext group.
Moreover, there is a factorization through the coinvariant ring
$[P_{\mathcal{A}} H_{*} (BV_{n})]_{GL(n)}$,
\[
\xymatrix{
 P_{\st } H_{*} (BV_{n})  \ar[rr]^{Tr_{n}^{\st}} \ar[dr]_{q}
                &  &   H^{n,n+*} (\st) ,  \\
                &  [P_{\st } H_{*} (BV_{n})]_{GL(n)} \ar@{.>}[ur]^{\varphi^{\st }_{n}}}
\]
and $\varphi^{\st} = \bigoplus_{n \geq 0} \varphi^{\st}_{n}$ is again a
homomorphism of a bigraded algebra. The map $\mathtt{Tr}^{\st}$ can be
thought of as the $E_{2}$ page of the stable transfer $B(\z2)^{n}_{+}
\to S^{0}$ (see Mitchell \cite{Mitchell-splitting}) hence the name
``transfer''. Singer computed this map in
small ranks, and found that $\varphi_{n}^{\st}$ is an isomorphism
for $n \leq 2$. Later, Boardman \cite{Boardman93}, with additional
calculations by Kameko \cite{Kameko}, showed that
$\varphi_{3}^{\st}$ is also an isomorphism. In fact, Singer has
conjectured that $\varphi^{\st}$ is always a monomorphism. Our main
interest in this paper is the image of $\varphi^{\st}$. It appears from the
calculations above that the image of the transfer $\mathtt{Tr}^{\st}$ is a
large, interesting and accessible subalgebra of $H^{*,*} (\st )$. In
particular, this image contains the Adams subalgebra on generators
$h_{i}$. Our calculation  strongly suggests  that $\varphi^{\st}$ also
detects many elements in the wedge subalgebra. In fact, in some sense,
elements in the wedge subalgebra have a better chance to be in the
image of the transfer than others.

The study of the Singer transfer is intimately related to the
problem of finding a minimal set of generators for the cohomology
ring $H^{*} (BV_{n})$ as a module over the Steenrod
algebra\footnote{This problem is called ``the hit problem'' in
the literature (see Wood \cite{Wood-problem}).}. The
$\st$--indecomposables in $H^{*} (BV_{n})$ were completely calculated
by Peterson~\cite{Peterson} for $n \leq 2$, and by
Kameko~\cite{Kameko} for $n=3$, and in ``generic degrees'' for all
$n$ by Nam~\cite{Nam-generic}. Here we prefer to work with the dual
$P_{\st} H_{*} (BV_{n})$ because of its ring structure, and also
because we are interested in the algebra structure by considering
this subring for all $n$.  We should mention that Smith and
Meyer~\cite{Smith-Meyer} have recently found a surprising connection
between the  subring $P_{\st} H_{*} (BV_{n})$  and a certain type
of  Poincar\'e duality quotient of the polynomial algebra, a
subject of great interest in modular invariant theory.

An important ingredient in Kameko's calculation~\cite{Kameko} of the
$\st$--generators for $H^{*}(BV_{3})$ is the existence of an operator
\[
Sq^{0}\co P_{\st} H_{d} (BV_{n}) \to P_{\st} H_{2d+n} (BV_{n}),
\]
for all $d , n \geq 0$. To explain the notation, recall that there
are Steenrod operations $\smash{\widetilde{Sq}^{i}}$ acting on the
cohomology of any cocommutative Hopf algebra (see May~\cite{May70}
or Liulevicius~\cite{Liulevicius60})  such that the operation
$\smash{\widetilde{Sq}^{0}}$ is not necessarily the identity. It turns
out that Kameko's operation commutes with $\smash{\widetilde{Sq}^{0}}$
via the Singer transfer (see Boardman~\cite{Boardman93}). This key
property is used by Bruner, H\`a and H\uhorn{}ng~\cite{BHH} to show that
the family of generators $g_{i} \in H^{4,*} (\st )$ is not in the
image of the transfer. As a result, $4$ is the first degree where
$\varphi_{4}^{\st}$ is not an epimorphism. In another direction,
Carlisle and Wood~\cite{Carlisle-Wood} showed that as a vector space,
the dimension of $P_{\st} H_{d} (BV_{n})$ is uniformly bounded; that is,
it has an upper bound which depends only on $n$. It follows that some
sufficiently large iteration of the endomorphism $Sq^{0}$ must become an
isomorphism. In fact, H\uhorn{}ng \cite{Hung-modular} showed that the number
of iterations needed is precisely $(n-2)$. This beautiful observation
allowed him to obtain many further information on the image of the Singer
transfer. Among other results, H\uhorn{}ng showed that for each $n \geq 5$,
$\varphi^{\st}_{n}$ is not an isomorphism in infinitely many degrees
(the same conclusion for $\varphi^{\st}_{4}$ can be deduced from the
main result of \cite{BHH}.) Moreover, using some computer calculations
by Bruner (using MAGMA)  and by Shpectorov (using GAP), H\uhorn{}ng also made
a comprehensive analysis of the image of the transfer in rank $4$  and
gave a conjectural list of elements in  $H^{4,*} (\st )$ that should,
or should not, be in the image of  $\varphi_{4}$ provided Singer's
conjecture is true in rank $4$.

Despite these successful calculations,  it seems that one has arrived
at the computation limit on both sides of the Singer transfer. What we
really need now are some global results on the structure of the graded
algebra $\bigoplus_{n} P_{\st} H_{*} (BV_{n})$.  This paper is the first
step in our investigation of the multiplicative structure of the graded
algebra $P_{\st} H_{*} (BV_{*})$ using suitably chosen sub-Hopf algebras
of $\st$. We shall show that the tranfer can be constructed not only
for the Steenrod algebra, but also any of its sub-Hopf algebras. We then
propose a conjecture which is the analog of the Quillen stratification
theorem for $P_{\st} H_{*} (BV_{*})$. We also make some calculations for
the transfer with respect to an important class of sub-Hopf algebras of
$\st$. One of the main results in this paper is the following application
to the study of the original Singer transfer that,
in our opinion, demonstrates the potential power of our approach.
\begin{theorem}\label{thm:1.1}
{\rm(i)}\qua The element $g \in H^{4, 24} (\st )$ is not in the
image of the Singer transfer.\\
{\rm(ii)}\qua The elements $d_{0} \in H^{4, 18} (\st)$ and $e_{0} \in
H^{4,21} (\st)$ are in the image of the Singer transfer.
\end{theorem}
We refer to Mahowald--Tangora~\cite{Mahowald-Tangora} and
Zachariou~\cite{Zachariou73} for detailed information about the generators
that appear in this
theorem. The fact that $g$, and in fact the whole family of generators $g_{i} =
\smash{(\widetilde{Sq}^{0})^{i}} g$, are not in the image of the Singer
transfer was already proven in \cite{BHH}. We give a different
proof, which is much less computational. The second part of our theorem
is new and give an affirmative answer to a part of H\uhorn{}ng's conjecture 
\cite[Conjecture 1.10]{Hung-modular}. It should be noted that our method seems
applicable to many other generators in the wedge subalgebra, but the
calculation seems much more daunting.

\subsubsection*{Organization of the paper}
The first two sections are preliminaries.
In \fullref{sec: intro}, we recall basic facts about the Steenrod algebra and its
sub-Hopf algebras. Detailed information about the action of the
operation $P_{t}^{s}$ on $H_{*} (BV_{n})$ is also given.  In
\fullref{sec: construction}, we present a convenient resolution, called the
\textit{Hopf bar resolution} by Anderson and Davis
\cite{Anderson-Davis}, to compute the cohomology of a Hopf
algebra. This resolution is then used to construct a chain-level representation of the Singer transfer for any sub-Hopf 
algebra B of $\st$. The idea of using this particular resolution to study the Singer transfer is due to Boardman \cite{Boardman93} (See also Minami \cite{Minami99}.) The remaining three sections are related, but independent of each other and can be read separately.  
We present our stratification conjectures for the domain of the Singer
transfer in \fullref{sec: Palmieri}.  
\fullref{sec: study of E transfer} is devoted to properties of the
$B$--transfer for various sub-Hopf algebras $B$ of $\st$.
\fullref{sec:info}
contains the proof of \fullref{thm:1.1}, which is one of the main applications of our approach. 

\subsubsection*{Acknowledgements}
I would like to thank Nguy\~\ecircumflex{}n H\,V H\uhorn{}ng for many
illuminating discussions and for sharing his ideas. I am also grateful
to Bob Bruner for many helpful conversations and to Nguy\~\ecircumflex{}n
Sum for pointing out an error in one of my calculations in earlier draft.

This article was partly conceived while I was a chercheur associ\'e
de CNRS at the Universit\'e de Paris Nord, and Universit\'e des Sciences et Technologies de
Lille. I would like to give special thanks to Daniel Tanr\'e and especially
Lionel Schwartz for their constant supports and encouragements.
\section{Sub-Hopf algebras of the Steenrod algebra}\label{sec: Milnor
generators} In this section, we review some basic facts about the
Milnor basis of the Steenrod algebra and the classification of sub-Hopf
algebras of $\st$. There are several excellent references on these
materials in literature; among them we highly recommend Margolis's book
\cite{Margolis-book} and Palmieri's memoir \cite{Palmieri-memoir}.
The original sources are Milnor \cite{Milnor} and Anderson and Davis
\cite{Anderson-Davis}.

\subsection{Milnor's generators}\label{subsec: Milnor}
It is generally more convenient to describe sub-Hopf algebras of the
Steenrod algebra in terms of their dual, as quotient algebras of
the dual Steenrod algebra $\st^{*}$. We recall some necessary
materials about the dual of the Steenrod algebra. According to
Milnor \cite{Milnor}, there is an algebra isomorphism
\[
\st^{*} \cong \f2 [\xi_{0}, \xi_{1}, \dots , \xi_n, \dots ],
\]
where $\xi_{t}$ is in degree $2^{t}-1$, and $\xi_{0}$ is understood to
be the unit $1$. The coproduct $\Delta$ is given by
\begin{equation}\label{eq: coproduct}
\Delta (\xi_n) = \sum_{i} \xi_{n-i}^{2^{i}} \otimes \xi_{i}.
\end{equation}
For any $s \geq 0$ and $t >0 $, let $P_{t}^{s} \in \st$ denote the dual of $\xi_{t}^{2^{s}}$. These
generators  are very important for our purpose. We briefly review some
of their fundamental properties.
If $s <t$, then $P_{t}^{s}$ is a differential, ie $(P^{s}_{t})^{2}
=0$. The effect of $P^{s}_{t}$ on $H^{*} (\mathbb{R}P^{\infty}) \cong
\f2[x]$ is completely determined by the formula
$$P^{s}_{t} x^{2^{k}}=x^{2^{s+t}}$$
if $k=s$, and zero for all $k \neq s$.
Let $b_{k} \in H_{*} (B\z2)$ denote the dual of $x^{k}$. We will be
working extensively with the dual action which reads
\begin{equation}\label{eq: action of P}
b_{k} P^{s}_{t} = \binom{k-2^{s} (2^{t}-1)}{2^{s}} b_{k - 2^{s} (2^{t}-1)},
\end{equation}
where binomial coefficients are taken modulo $2$. Write $2^{s} \in k$ if
$2^{s}$ appears in the binary expansion of $k$ and $2^{s} \notin k$ if
the opposite happens. We will need the following simple but useful lemma.
\begin{lemma}\label{lem: Pst action} With the notation as above,
\begin{itemize}
\item [(i)] $b_{k} P_{t}^{s}=0$ if and only if either $k <
2^{s+t}$, or $k \geq 2^{s+t}$ and $2^{s} \in k$.
\item [(ii)] $b_{k}$ is in the image of $P_{t}^{s}$ if and only if $k \geq
2^{s+t}$ and $2^{s} \in k$.
\end{itemize}
\end{lemma}
\begin{proof}
$\smash{\tbinom{k-2^{s} (2^{t}-1)}{2^{s}}}$ is non-zero if and only if $k
\geq 2^{s+t}$, and $2^{s} \in k-2^{s} (2^{t}-1)$. The latter
condition is clearly equivalent to $2^{s} \notin k$.
\end{proof}


\subsection{Sub-Hopf algebras of the Steenrod algebra}\label{sec:
subHopf}
We are mainly interested in two families of sub-Hopf algebra of $\st$:
the elementary ones, which essentially play similar role as
elementary abelian subgroups in group cohomology; and the normal ones,
which serve as intermediate between the elementary sub-Hopf algebra and the whole
Steenrod algebra.

Let $A$ be a Hopf algebra. A sub-Hopf algebra $E \subset A$ is called 
\emph{elementary} if it is bicommutative, and $e^{2} =0$ for any element e in the augmentation
ideal $\bar{E}$ of $E$. This definition is due to Wilkerson
\cite[page~138]{Wilkerson-TAMS81}. We now specialize to the case $A = \st$.  Each
elementary sub-Hopf algebra $E$ of $\st$ is  isomorphic, as an
algebra, to the exterior algebra on the operation $P^{s}_{t}$s that  it contains. In particular,
\begin{equation}\label{eq: coho E}
H^{*} (E) \cong \f2 [h_{ts}| P^{s}_{t} \in E],
\end{equation}
where $h_{ts}$ is represented by $[\xi_{t}^{2^{s}}]$ in the cobar
complex for $E$, so $|h_{ts}|= (1, 2^{s}(2^{t}-1))$. Among these
elementary sub-Hopf algebras, the maximal ones have the form
\begin{equation}\label{eq: dual of E(n)}
E(m)^{*} = \st^{*}/(\xi_1 , \dotsc ,
\xi_{m-1}, \xi_{m}^{2^{m}}, \xi_{m+1}^{2^{m}}, \dotsc ),
\end{equation}
for each $m \geq 1$. Equivalently, $E(m)$ is generated by those $P_{t}^{s}$ for which
$s < m \leq t$.
We now discuss normal sub-Hopf algebras of $\st$.
We say that a sub-Hopf algebra $B$ is normal in $\st$ if the left and
the right ideal generated by $\bar{B}$ are equal; that is,  $\bar{B}
\st = \st \bar{B}$.  If
$B$ is normal in $\st$, then one can define the quotient Hopf algebra
$\st {\dslash} B$ as $\st \otimes_{B} \; \f2 = \f2 \otimes_{B} \st $. The short
exact sequence of vector spaces $B \to \st \to \st {\dslash} B$ is called a
\textit{Hopf algebra extension}. Of course, this definition applies
not just for the Steenrod algebra, but also to any cocommutative Hopf algebra over
$\f2$. The normal sub-Hopf algebras of $\st$ are completely
classified (see Margolis \cite[Theorem 15.6]{Margolis-book}).
They correspond to non-decreasing sequences $n_{1} \leq n_{2} \leq \dotsb
\leq \infty$ via the correspondence
\[
(n_{1}, n_{2}, n_{3}, \dotsc) \to  \st^{*}/(\xi_{1}^{2^{n_{1}}},
\xi_{2}^{2^{n_{2}}}, \xi_{3}^{2^{n_{3}}}
\dotsc ).
\]
In particular, maximal elementary sub-Hopf algebras
$E(m)$ are normal in $\st$. The union of $E(m)$s, denoted by $D$, turns out to be another
normal sub-Hopf algebra of $\st$. In fact, as observed by
Palmieri~\cite{Palmieri99},  there is a whole sequence
of normal sub-Hopf algebras, starting at $D$
\begin{equation}\label{eq: D sequence}
D= \bigcap_{m} D(m) \to \dotsb  \to D(m) \to D(m-1) \dotsb  \to D(1) \to D(0) = \st,
\end{equation}
where $D(m)$ is defined in terms of its dual as follows,
\begin{equation}\label{eq: definition of D}
D(m)^{*}= \st^{*}/(\xi_1^{2}, \xi_2^{4}, \dotsc , \xi_m^{2^{m}}).
\end{equation}
In other words, $D(m)$ is generated by the operations $P^{s}_{t}$s
where either $t > m$ or $s < t \leq m$. The quotient $D(m-1){\dslash}D(m)$ is
the exterior algebra on generators $P_{m}^{m+i}$, $i \geq 0$. In
particular, $D(m-1){\dslash}D(m)$ is ($2^{m} (2^{m}-1) -1$)-connected.

\section{Construction of the $B$--transfer}\label{sec: construction}
Let B be any sub-Hopf algebra of $\st$. Clearly B is also a graded
connected and cocommutative, so that results from  \fullref{sec: construction}
can be applied to B. We will construct a chain level representation $P_{B} H_{*} (BV_{n}) \to (\bar{B}^{*})^{n}$ for 
the analogue of the Singer transfer for any sub-Hopf algebra $B$ of
$\st$. We begin with a review of the so-called (reduced) Hopf bar resolution that we will use. 
\subsection{Hopf bar resolution}\label{subsec: Hopf resolution}
Let $A$ be a graded cocommutative connected Hopf
algebra over $\f2$ and let $M$ be an $A$--module. We denote by $\mu \co A
\otimes A \to A$ the product and $\Delta \co A \to  A \otimes A$ the
coproduct maps. In this section, we
present the (normalized) Hopf bar construction of $A$ for $M$, introduced
by Anderson and Davis~\cite{Anderson-Davis},  to calculate the cohomology $H^{*,*}(A, M)$ of
$A$  with coefficient in $M$. This construction is functorial with respect to maps between Hopf
algebras as well as maps between $A$--modules when $A$ is fixed. This particular
resolution is well-suited to our purposes rather than the usual bar
resolution because, as we shall see later, there exists an explicit description of the
representing map for the transfer. 

Observe that if $A$ is a cocommutative connected Hopf algebra, then the
tensor product of any two $A$--modules is again an $A$--module via the
coproduct $\Delta$.  Let $\bar{A}$ be the augmentation ideal of $A$. From the obvious
short exact sequence of $A$--modules $0 \rightarrow \bar{A}
\rightarrow A \to \f2 \to 0$, tensoring (over $\f2$) with
$\bar{A}^{\otimes k} \otimes M$ (from now on, we will write $\bar{A}^{k}$ instead
of $\bar{A}^{\otimes k}$ to avoid clustering) and splicing together the
resulting short exact sequences, we obtain a chain complex
$\mathcal{H} (M)$ which is visibly exact:
\begin{equation}\label{eq: exact complex}
\dots \to A \otimes \bar{A}^{k} \otimes M \to A \otimes
\bar{A}^{(k-1)} \otimes M \to \dots \to A \otimes M \to M.
\end{equation}
We claim that $\mathcal{H}(M)$ is an $A$--free resolution of $M$. Indeed, it
suffices to verify that $A \otimes \bar{A}^{k} \otimes
M$ is a free $A$--module for each $k \geq 0$. This is not quite as obvious
as it seems because by its construction, the $A$--module structure of $A \otimes \bar{A}^{k} \otimes
M$  is via the iterated coproduct. However, this $A$--action can be
modified by mean of a well-known trick that for any $A$--module $N$, the
usual (ie $A$--module structure via coproduct) $A$--module  $A
\otimes N$ is isomorphic to the $A$--module $A \otimes tN$ where $A$ acts
only on the copy of $A$, and $tN$ signifies the same $\f2$--vector space
$N$, but with trivial $A$--action (see Anderson and Davis
\cite[Proposition~2.1]{Anderson-Davis}).

Let $M^{*} = \Hom_{\f2} (M, \f2)$ be the $\f2$--linear dual of $M$. Taking $\Hom_{A}
(\mathcal{H}(M), \f2)$  and simplify, we obtain a cochain complex 
%
\[
0 \rmap{} M^{*} \rmap{\partial_{0}} \bar{A}^{*} \otimes M^{*}
\rmap{\partial_{1}} \cdots \rmap{}
(\bar{A}^{*})^{k}\otimes M^{*}  \rmap{\partial_{k}}
(\bar{A}^{*})^{(k+1)}\otimes M^{*}  \rmap{} \cdots.
\]
whose homology is  $H^{*,*} (A; M)$. Since we have modified the
$A$--action on $A \otimes \bar{A}^{k} \otimes M$,  the differential
$\partial_0$ becomes a kind of twisted coaction map,
\begin{align}\label{eq: diff}
\partial_{0}\co M^{*} \rmap{\alpha_*} \bar{A}^{*} \otimes M^{*} \rmap{\chi \otimes id}
\bar{A}^{*} \otimes M^{*},
\end{align}
where $\alpha\co A \otimes M \to M$ is the map that defines the $A$--module structure on
M.
There is similar description for $\partial_{k}$, $k
>0$. Namely, if $m \in M^{*}$ with $ \alpha_{*}
(m)= \sum_{\nu } a_{\nu } \otimes m_{\nu }$, then
\begin{equation}\label{eq: k differential}
\partial_{k} (a_{1}|\dotsb |a_{k}|m) = \sum_{\nu} \sum \chi (a_{1}'\dotsb
a_{k}' a_{\nu })|a_{1}''|\dotsb |a_{k}''|m_{\nu }.
\end{equation}
if $\mu^{*} (a_{i}) = \sum a_{i}' \otimes a_{i}''$ (See Anderson and
Davis~\cite[pages 320--321]{Anderson-Davis}.)
In particular, if $M=\f2$, then equation~\eqref{eq: k differential}
becomes 
\begin{equation}\label{eq: diff when M=F2}
\partial_{k} (a_{1}|\dotsb |a_{k}) = \sum (\chi (a_{1}'\dotsb
a_{k}')|a_{1}''|\dotsb |a_{k}'').
\end{equation}

\subsection{Construction of the $B$--transfer}\label{subsec:
construction}
We begin with the existence of the $B$--transfer.
\begin{theorem}\label{thm: Singer main result} Let $B$ be a sub-Hopf
algebra of $\st$. For each $n \geq 0$, there exists a map
\[
P_{B} H_{*} (BV_{n}) \rmap{\mathtt{Tr}_{n}^{B}} H^{n,n+*} (B),
\]
natural with respect to $B$, such that it factors through the
coinvariant ring
\[
\xymatrix{
 P_{B} H_{*} (BV_{n})  \ar[rr]^{\mathtt{Tr}_{n}^{B}} \ar[dr]_{q}
                &  &    H^{n,n+*} (B).    \\
                &  [P_{B} H_{*} (BV_{n})]_{GL(n)} \ar@{.>}[ur]^{\varphi^{B}_{n}}}
\]
Furthermore, the total $B$--transfer $\varphi^{B} = \bigoplus_{n} \varphi^{B}_{n}$
is a homomorphism of bigraded algebras.
\end{theorem}
\begin{proof}
A careful look at Singer's construction shows that his proof also
works for any sub-Hopf algebra $B$ of $\st$. We follow
Boardman's approach because it provides an explicit description of a map
$P_{B} H_{*} (BV_{n}) \to  (\bar{B}^{*})^{n}$ which represents
$\mathtt{Tr}_{n}^{B}$. 

To begin, we need to introduce several $\st$--modules related to $P =
\f2[x] = H^{*} (B\z2)$. Let $\hat{P}$ be obtained from
$P$ by formally adding a basis element $x^{-1}$ in degree $-1$ and equip
$\hat{P}$ with an $\st$--action such that it is the unique extension of
the $\st$--action on $P$ that also satisfies the Cartan
formula. According to the referee, $\hat{P}$ was
first introduced by Adams in \cite{Adams-1972}. 

Let $\hat{f}$ be the $\st$--epimorphism $\st \to \hat{P}$ such that $\hat{f}(1) = x^{-1}$. Denote by $f$  
its restriction to $\wwbar{\st}$.  It is clear that $f$ maps into
$P$. Moreover, $f$ is $B$--linear for any sub-Hopf algebra $B$ of
$\st$. On the other hand, the inclusion $B \to \st$ provides $\st$
with a $B$--module structure,  so the dual $\sst \to B^{*}$ is a map of
right $B$--modules.
It follows that the following composition is a map of right
$B$--modules:
\begin{equation}\label{eq: chain level map}
f_{n}^{B} \co H_{*} (BV_{n}) \rmap{f^{\otimes n}_{*}}
(\wwbar{\st}^{*})^{n} \to (\bar{B}^{*})^{n}.
\end{equation}
The image of this map of any $B$--annihilated element in $H_{*}
(BV_{n})$ is a cocycle in the Hopf resolution for $B$. Thus
$f_{n}^{B}$ induces a map
\[
\mathtt{Tr}_{n}^{B}\co P_{B} H_{*} (BV_{n}) \rmap{} H^{n,n+*} (B),
\]
which is our version of the Singer transfer for the sub-Hopf algebra $B$. It is
clear that the construction just described is natural with respect to $B$.
\end{proof}

Singer's proof that $\mathtt{Tr}_{n}$ factors through $GL(n)$
coinvariants is very simple and elegant. Recall that $GL(n)$ is generated
by the symmetric group $\Sigma_{n}$ and an element denoted by $\tau$
of order $3$ in $GL(2)$, considered as a subgroup of $GL(n)$ in the
obvious way. In terms of the chain-level map that we have just
constructed, the fact that $\mathtt{Tr}_{n}$  factors through
$\Sigma_{n}$--coinvariants is precisely because the  Steenrod algebra
$\st$ is cocommutative.  That $\mathtt{Tr}_{2}$ is invariant under
$\tau$ seems much less obvious.  In fact, we have to use a smaller
resolution, which is the Lambda algebra. The author plans to write
about this elsewhere. 

\begin{remark}\label{rem: image of f} According to Boardman
\cite{Boardman93}, the $\st$--linear map $f_{*}$ has an explicit
description as follows. Let $b_{k}$ be the generator of
$H_{k} (B\z2)$ in degree $k$, dual to $x^{k} \in H^{k} (B\z2)$. Then
the image of $b_{k}$ under $f_{*}$ is the coefficient of $x^{k+1}$ in the expansion
of
\begin{equation}\label{eq: Boardman formula}
\prod_{i=0}^{\infty} (1 \otimes 1 + x^{2^{i} (2^{1}-1)} \otimes
\xi_{1}^{2^{i}} + x^{2^{i}(2^{2}-1)} \otimes
\xi_{2}^{2^{i}} + x^{2^{i}(2^{3}-1)} \otimes
\xi_{3}^{2^{i}}  + \dotsb ).
\end{equation}
\end{remark}

Of course, for a fixed degree $k$, one needs only to consider finite
products. The situation is even simpler for $f_{n}^{B}$ when $B$ is
small because many elements of the form $\xi_{t}^{2^{s}}$ get killed when mapped
down to $\bar{B}^{*}$.

For example, if $B= E(1)$, then the only nontrivial factor in the
above infinite product is $i=0$ and so the only nontrivial images are $f_{*}(b_{2^{t}-2}) = \xi_{t}$, $t
\geq 1$.

Now let $\mathcal{E}$ be the collection of all elementary sub-Hopf
algebras of $\st$. By the naturality of the transfer $\mathtt{Tr}^{B}$,
we have the commutative diagram
\begin{equation}\label{eq: comm diag 1}
\xymatrix{
 P_{\st} H_{*} (BV_{n}) \ar[d]_{\mathtt{Tr}^{\st}_{n}} \ar[r]^{i^{\st}_{D}}
                & P_{D} H_{*} (BV_{n}) \ar[d]_{\mathtt{Tr}^{D}_{n}}
\ar[r]^<<<<{i^{D}_{\mathcal{E}}}_<<<{\cong} &  \invlim{\mathcal{E}} P_{E} H_{*}
(BV_{n})\ar[d]^{\lim \mathtt{Tr}^{E}_{n}}  \\
H^{n,n+*} (\st ) \ar[r]^{\res^{\st}_{D}} & H^{n,n+*} (D )\ar[r]^<<<<{\res^{D}_{\mathcal{E}}} & \invlim{\mathcal{E}}  H^{n,n+*} (E ),       }
\end{equation}
where the horizontal maps are, or are induced by, the obvious inclusions. Note
that $i^{\st}_{D}$ is a monomorphism, and $i^{D}_{\mathcal{E}}$ is
an isomorphism because $D$ is generated by the elementary sub-Hopf
algebras of $\st$. Of course, the induced maps  $\varphi_{n}^{B}$ after passing to
$GL(n)$--coinvariant rings  for various $B$s in the diagram above between
need not be mono nor epi.

For convenience, write $\mathtt{Tr}^{\mathcal{E}}_{n}$
for the inverse limit $\smash{\invlim{E}} P_{E} H_{*} (BV_{n})$, and call it the
$\mathcal{E}$--transfer. The factorization of the
$\mathcal{E}$--transfer to the coinvariant ring is denoted by
$\varphi^{\mathcal{E}}$. 

\subsection{Kameko's $Sq^{0}$}\label{subsec: sqo on cohomology}
Perhaps the the most useful tool in the study of the hit problem and the Singer transfer is a certain operator defined by Kameko \cite{Kameko}.  We discuss the behavior of this operator with respect to the $B$--transfer.   
According to Liulevicius \cite{Liulevicius60} (see also May
\cite{May70}), the cohomology of any cocommutative Hopf algebra
$A$ is equipped with an action of the Steenrod algebra where
$\smash{\widetilde{Sq}^{0}}$ may act nontrivially (ie it is not necessarily the
identity). In fact, $\smash{\widetilde{Sq}^{0}} \co \smash{H^{n,q}} (A) \to
\smash{H^{n,2q}} (A)$ is induced by the Frobenius $z \mapsto z^{2}$ in the
cochain level.

Kameko's operation on $P_{\st} H_{*} (BV_{n})$ behaves much like
$\smash{\widetilde{Sq}^{0}}$. By definition, it is a map $H_{d} (BV_{n}) \to
H_{2d+n} (BV_{n})$, given by the formula
\[
Sq^{0} \co b_{i_{1}} \dotsb b_{i_{n}} \mapsto
b_{(2i_{1} +1)} \dotsb b_{(2i_{n}+1)}.
\]
In fact, it is easy to check that under the inclusion $f_{*}$, this
operation corresponds precisely with the Frobenius homomorphism in
$(\sst)^{n}$.
The following easy lemma describes the relation between Kameko's $Sq^{0}$ and the
operations $P_{t}^{s}$.
\begin{lemma}\label{lem: Pst and Sq0}  $(Sq^{0}z) P_{t}^{0} = 0$ and  $(Sq^{0}z)P^{s}_{t} = Sq^{0}
(zP^{s-1}_{t})$ when $s > 1$ for any $z \in H_{*}
(BV_{n})$.
\end{lemma}
\begin{proof} This is immediate from formula~\eqref{eq: action of P}.
\end{proof}
It follows from this lemma that Kameko's operation also induces an
automorphism on $P_{B} H_{*} (BV_{n})$ where B is one of the sub-Hopf
algebras $E(m), D(m)$ or D. Thus Kameko's $Sq^{0}$ commutes with the
Liulevicius--May $\smash{\widetilde{Sq}^{0}}$ via the $B$--transfer for any
sub-Hopf algebras of the types above.

\begin{remark}
In  \cite{Carlisle-Wood}, Carlisle and Wood proved a striking
property that the dimension of the vector space $P_{\st} H_{*}
(BV_{n})$ is bounded,  for each $n$ fixed.  It can be shown that the
same is true when replacing $\st$ by any sub-Hopf algebra $D(m)$.
This is essentially because the quotient  algebra $\st{\dslash}D(m)$ is
finite. On the other hand, $P_{D} H_{*} (BV_{n})$ is not uniformly
bounded.  Here is an example in rank 2. Choose a number $a >1$, and
for all $1 \leq i \leq a-1$, let $k_{i} =   2^{a} (2 \times 2^{i}
-1) -1$, and $\ell_{i} =  2^{a+i+1} (2 \times 2^{a-i -1} -1) -1$.
One can verify, using \fullref{cor: D-annihilated} below, that
in degree $d= 2^{2a+1} -2^{a} -1$, all $(a-1)$ monomials $b_{k_{i}}
b_{\ell_{i}}$  are $D$--annihilated. Since $a$ is chosen arbitrarily,
it follows that $P_{D} H_{*} (BV_{2})$ is not uniformly bounded.
\end{remark}

\section{Two stratification conjectures}\label{sec: Palmieri}
In this section, we discuss our conjectures on the domain of the Singer
transfer, which are the analogues of Quillen stratification theorem about
the cohomology of finite groups, and Palmieri's version for cohomology
of the Steenrod algebra. These conjectures, if true, would provide some
global information on the mysterious algebra $\bigoplus_n P_{\st} H_{*}
(BV_n)$.

The lower horizontal maps in ~\eqref{eq: comm diag 1} were studied by Palmieri
\cite{Palmieri99}. We shall first summarize his results.

\subsection{Palmieri's stratification theorems}
Just as in the classical case of group cohomology, there is an action
of the Hopf algebra $\st$ on the cohomology of its
sub-Hopf algebra $H^{*,*} (D)$ such that the image of the restriction
map $\res^{\st}_{D}$ actually lands in the subring $[H^{*,*}
(D)]^{\st}$ of elements that are invariant under this action. Since
$D$ is normal in $\st$ and since $D \subset \st$ acts trivially on the
cohomology of itself, we can write $H^{*,*}
(D)^{\st} = H^{*,*} (D)^{\st{\dslash}D}$. The following is
Palmieri's version of Quillen stratification for the Steenrod algebra.
\begin{theorem}[Palmieri~\cite{Palmieri99}]\label{thm: P theorem I}
The canonical maps
\[
H^{*,*} (\st) \rmap{\res^{\st}_{D}} [H^{*,*} (D)]^{\st{\dslash}D},
\]
and
\[
H^{*,*} (D) \rmap{\res^{D}_{\mathcal{E}}} {\textstyle\invlim{\mathcal{E}}}  H^{*,*} (E ),
\]
are $F$--isomorphisms. That is, their kernel and cokernel are nilpotent
(in an algebraic sense).
\end{theorem}
Moreover, the limit $R= \smash{\invlim{\mathcal{E}}}  H^{*,*} (E )$ is
computable. To describe it, recall that the cohomology of each
elementary sub-Hopf algebra $E$ is a polynomial algebra on  elements
$h_{t,s}$ which corresponds to each generator $P_{t}^{s}$ that it
contains.  The second theorem of Palmieri says:
\begin{theorem}[Palmieri {\cite[Theorem 1.3]{Palmieri99}}]
\label{thm: P inverse limit} 
There is an isomorphism of $\st$--algebras
\[
R = \f2 [h_{t,s}|s <t]/(h_{t,s}
h_{v,u}|u \geq t),
\]
where $|h_{t,s}|= (1, 2^{s}(2^{t}-1))$. The action of $\st$ is given by
the Cartan formula and the following formula on the generators
\[
Sq^{2^{k}} h_{t,s}=
\begin{cases}
h_{t-1,s+1} &   \text{if $k=s$ and $t-1 >s+1$},\\
h_{t-1,s} & \text{if $k=s+t-1$ and $t-1 > s$},\\
0 & \text{otherwise.}
\end{cases}
\]
\end{theorem}
It is important to note that the $\st$--action here is induced by the inclusion $E
\subset \st$  for each $E \in \mathcal{E}$, not the kind of
$\st$--action of Liulevicius--May on cohomology of cocommutative Hopf
algebras that we mentioned earlier. For the latter action, it can be
verified in the cochain level that $\smash{\widetilde{Sq}^{0}}$ acts
on $R$ by sending $h_{t,s}$ to $h_{t,s+1}$ if $s+1 <t$, and sending
$h_{t,t-1}$ to zero.

Combining the two F-isomorphisms in \fullref{thm: P theorem I}, we
obtain an isomorphism mod nilpotents $H^{*,*} (\st) \to R^{\st}$. There are two major
difficulties if we want to use this map to study $H^{*,*}
(\st)$. First of all, the invariant ring $R^{\st}$, which is another kind of "hit problem", seems to be very
complicated because it has lots of zero divisors. Secondly, given an invariant element in $R^{\st}$, we do not
know which power of it is lifted to $H^{*,*} (\st)$. On the other hand,
several families of elements in $R^{\st}$ are known and we will discuss
the problem whether the Singer transfer or its $B$--analogues detects several
such families in Sections \ref{sec: study of E transfer} and \ref{sec:info}.

\subsection{Two conjectures}\label{subsec: conjectures}
Base on Palmieri's results and the Singer conjecture, we make the
following conjectures.

\begin{conjecture}\label{conj: my own}
The following canonical homomorphism of algebras,
\begin{align}\label{eq: conj}
i^{\st}_{D}\co \bigoplus_{n} P_{\st} H_{*} (BV_{n})_{GL(n)} \rmap{} &
\bigoplus_{n}
[P_{D} H_{*} (BV_{n})_{GL(n)}]^{\st{\dslash}D}, \quad  \text{and}\\
i^{D}_{\mathcal{E}}\co \bigoplus_{n} [P_{D} H_{*} (BV_{n})]_{GL(n)}
\rmap{} & \bigoplus_{n} {\textstyle\invlim{E \in
\mathcal{E}}} [P_{E} H_{*} (BV_{n})]_{GL(n)},
\end{align}
are both $F$--isomorphisms.
\end{conjecture}
Here is our evidence for making the above conjecture. First of all,
they are true for trivial reasons in degree $n=1$. In fact, both maps
are isomorphisms in this degree. Secondly, recall the Singer conjecture
that $\varphi^{\st}$ is a monomorphism. If this conjecture is true,
then according to \fullref{thm: P theorem I}, it is necessary that
the kernel of the two homomorphisms above be nilpotent. Our final
evidence is the fact that $\varphi_{n}^{\st}$ is an isomorphism for $n
\leq 3$. But in these degrees, only $h_{0} \in H^{1,1} (\st)$ has
infinite height, so the kernel of $i^{\st}_{\mathcal{E}}$ is
nilpotent for $n \leq 3$.

We now describe a refinement, by adapting Palmieri's proof of
his theorem, of the first part of the above conjecture. If $M$ is an
$\st$--module of finite type, then the family of normal sub-Hopf
algebras $D(m)$ induces a filtration
\[
P_{D} M \supseteq \dotsb \supseteq P_{D(m)} M \supseteq P_{D(m-1)} M
\supseteq \dotsb \supseteq P_{\st} M,
\]
which stabilizes in each degree. The first part of \fullref{conj: my own}
will be a consequence of the following.
\begin{conjecture}\label{conj: first part modified} For each $m\geq
1$. The canonical map
\[
\bigoplus_{n} [P_{D(m-1)} H_{*} (BV_{n})]_{GL(n)} \to \bigoplus_{n}
([P_{D(m)} H_{*} (BV_{n})]_{GL(n)})^{D(m-1){\dslash}D(m)},
\]
is an $F$--isomorphism.
\end{conjecture}
Of course, $P_{D(m-1)} H_{*} (BV_{n}) = [P_{D(m)} H_{*}
(BV_{n})]^{D(m-1){\dslash}D(m)}$. The problem alluded to is the order of taking
$GL(n)$--coinvariant
and taking invariant under Steenrod operations.

For the second part of \fullref{conj: my own}, we have the
following observation. Since the family of maximal elementary sub-Hopf
algebra $E(m)$ is cofinal in the category $\mathcal{E}$, we can define
an element in $\smash{\invlim{E}} [P_{E} H_{*} (BV_{n})]_{GL(n)}$ as
a family of compatible elements $z_{m} \in [P_{E(m)} H_{*}
(BV_{n})]_{GL(n)}$. Compatibility means that the restriction of $z_{m}$
and $z_{m'}$ to $[P_{E(m) \cap E(m')} H_{*} (BV_{n})]_{GL(n)}$ must be
the same for any $m$ and $m'$.

On the other hand, in each fixed degree $d$, any element of $H_{d}
(BV_{n})$ is $E(m)$--annihilated for $m$ large. It follows that in a
compatible family $\{z_{m} \}$, one has
\[
z_{m} = z_{m+1} = \dotsb \in [H_{d} (BV_{n})]_{GL(n)},
\]
for $m$ sufficiently large. It is well-known that the algebra $[H_{*}
(BV_{n})]_{GL(n)}$, which is dual to the Dickson algebra of
$GL(n)$--invariants of polynomial algebra on $n$ generators, is
$(2^{n-1}-2)$--connected. It follows that any element in the algebra
$$\bigoplus_{n} [H_{*} (BV_{n})]_{GL(n)}$$
is nilpotent. Indeed, if $z \in
[H_{d} (BV_{n})]_{GL(n)}$, then $z^{k} \in [H_{dk}
(BV_{nk})]_{GL(nk)}$ is in degree $dk$, which is less that the
connectivity $(2^{nk-1}-2)$ for $k$ large.

We have shown that if
$$\bigoplus_{n}{\textstyle\lim_{E}^{f}}[P_{E}H_{*}(BV_{n})]_{GL(n)}$$
denotes the subalgebra of
$$\bigoplus_{n}{\textstyle\invlim{E}}[P_{E}H_{*}(BV_{n})]_{GL(n)}$$
consisting of \textit{finite} sequences $\{z_{m} \}$ (ie $z_{m} =0$
for all but finitely many values of $m$). Then
\begin{lemma}\label{lem: refinement of second part}
The inclusion of algebras
\[
\bigoplus_{n} {\textstyle\invlim{E}^{f}}[P_{E}H_{*}(BV_{n})]_{GL(n)}
\to \bigoplus_{n}
{\textstyle\invlim{E}}[P_{E}H_{*}(BV_{n})]_{GL(n)},
\]
is an $F$--isomorphism.
\end{lemma}
This result is relevant to the first part of the conjecture above as
well, for there is a theorem of Hung and Nam \cite{Hung-Nam,Hung-Nam-JA} which
says that the canonical homomorphism
\[
[P_{\st} H_{*} (BV_{n})]_{GL(n)} \to P_{\st} ([H_{*} (BV_{n})]_{GL(n)}),
\]
is trivial in positive degrees, as soon as $n \geq 3$. Thus the image
of the canonical map
\[
\bigoplus_{n} [P_{\st} H_{*} (BV_{n})]_{GL(n)} \to \bigoplus_{n}
{\textstyle\invlim{E}} [P_{E} H_{*} (B(BV_{n}))]_{GL(n)},
\]
actually lands in the subalgebra of finite sequences in
\fullref{lem: refinement of second part}.

It would be interesting also to see whether the theorem of Hung and
Nam mentioned above remain true when the Steenrod algebra is replaced
by a sub-Hopf algebra, say the family $D(m)$, or even the sub-Hopf
algebra $D$.

\section{Study of the $E$--transfer}\label{sec: study of E transfer}
In this section, we investigate the subring $P_{E} H_{*} (BV_{n})$
where $E$ is an elementary sub-Hopf algebra of $\st$. We pay
particular attention to the case when $E$ is a maximal one.
\subsection{$E(m)$--annihilated elements in $H_{*}
(B\z2)$}\label{subsec: rank 1 case}
The formula \eqref{eq: Boardman formula} serves as an efficient means
to compute the image of $b_{k}$ in $\st^{*}$ as well as
$E^{*}$. Recall that we have $f_{*} (b_{2k+1}) =
(f_{*} b_{k})^{2}$, so in principle, we need only to compute
the image of evenly graded elements.

\begin{notation}\label{not: kappa et rho} For each $k \geq 0$, write
$k +1 = 2^{\kappa} (2\rho - 1)$. Clearly, $\kappa$ and $\rho$ are
uniquely determined by $k$. In fact, $\kappa$ is the
smallest non-negative integer such that $2^{\kappa}
\notin k$. We reserve the letters $\kappa $ and $\rho $ for such a
presentation.
\end{notation}
\begin{lemma}\label{lem: E(m)-annihilated}
$b_{k}$ is $E(m)$--annihilated if and only if either (i) $\kappa \geq
m$ or (ii)  $\kappa < m$ and $\rho  \leq 2^{m-1}$.
\end{lemma}
\begin{proof} Recall that $E(m)$ is generated by the operations $P_{t}^{s}$
in which $s < m \leq t$. Thus to be $E(m)$--annihilated, $k$ must
satisfies the condition that for any $s < m$, if $2^{s} \notin k$,
then $k < 2^{s+t}$ for all $t \geq m$.

If $\kappa  \geq m$, then $2^{s} \in k$ for any $s < m$, so $b_{k}$ is
clearly $E(m)$--annihilated. If $\kappa < m$, then $k = 2^{\kappa} (2\rho-1) -1
< 2^{\kappa +m}$ which implies $\rho  \leq  2^{m-1}$.
\end{proof}

The following Corollary is immediate.
\begin{corollary}\label{cor: D-annihilated} $b_{k}$ is $D$--annihilated
if and only if $\rho  \leq  2^{\kappa}$.
\end{corollary}

\subsection{Some properties of the $E(m)$--transfer}\label{subsec: rank
n}
Recall from the proof of \fullref{thm: Singer main result} that
the composition
\[
f_{n}^{E(m)} \co P_{E(m)} H_{*} (BV_{n}) \rmap{f_{*}^{\otimes
n}} (\wwbar{\st}^{*})^{n} \to (\overline{E(m)}^{*})^{n},
\]
is a chain level representation for $\mathtt{Tr}^{E(m)}_{n}$. Moreover,
\[
E(m)^{*} = \st^{*}/(\xi_{1}, \dotsc , \xi_{m-1}, \xi_{m}^{2^{m}},
\xi_{m+1}^{2^{m}}, \dotsc ).
\]
Our first result says that the image of $f_{n}^{E(m)}$ has a rather
strict form.
\begin{lemma}\label{lem: image of Tnm} Under $f_{1}^{E(m)}\co H_{*}
(B\z2) \to \overline{E(m)}^{*}$, the image of $b_{k}$
is nontrivial if and only if $k$ can be written in the form
\begin{equation}\label{eq: nontrivial image}
k = 2^{s_{1}} (2^{t_{1}}-1) + \dotsb  +  2^{s_{\ell}}(2^{t_{\ell}}-1)-1,
\end{equation}
for some $0 \leq s_{1} < \dotsb < s_{\ell} \leq (m-1)$, and $m \leq t_{1}, \dotsc ,
t_{\ell}$. Moreover, the $\ell$--tuple $(s_{1},\dotsc s_{\ell})$ is
unique, up to a permutation.
\end{lemma}
\begin{proof}
We make use of Boardman's formula \eqref{eq: Boardman formula}. When
projecting down to $\overline{E(m)}^{*}$, the infinite product in this
formula is actually a finite product since $\xi_{t}^{2^{s}} = 0$ in
$\overline{E(m)}_{*}$ for all $s \leq m$. So the image of $b_{k}$ under
$T_{1}^{E(m)}$ is the coefficient of $x^{k+1}$ in the finite product
\begin{equation*}
\prod_{i=0}^{m-1} (1 \otimes 1 + x^{2^{i} (2^{m}-1)} \otimes
\xi_{m}^{2^{i}} + x^{2^{i}(2^{m+1}-1)} \otimes
\xi_{m+1}^{2^{i}} + x^{2^{i}(2^{m+2}-1)} \otimes
\xi_{m+2}^{2^{i}}  + \dotsb ).
\end{equation*}
Thus $b_{k}$ is mapped to the sum of $\prod_{i=0}^{\ell} \xi_{t_{i}}^{2^{s_{i}}}$
for each presentation of $k$ in the form \eqref{eq: nontrivial image}.
For the uniqueness of the set $\{s_{1},\dotsc ,s_{\ell} \}$. Observe
that the equality
\[
 2^{s_{1}} (2^{t_{1}}-1) + \dotsb  +  2^{s_{i}} (2^{t_{i}}-1) =
 2^{u_{1}} (2^{v_{1}}-1) + \dotsb  +  2^{u_{j}} (2^{v_{j}}-1),
\]
implies that
\[
2^{s_{1}} + \dotsb + 2^{s_{i}} \equiv 2^{u_{1}} + \dotsb +
2^{u_{j}} \pmod 2^{m},
\]
But both sides in the above equation are at most
$2^{0} + \dotsb  + 2^{m-1} = 2^{m} -1$, hence they must be
actually equals. We then obtain two binary expansion of
the same number, which implies that the two set of indices
$\{s_{1},\dotsc ,s_{i} \}$ and $\{u_{1},\dotsc ,u_{j} \}$ are the
same.
\end{proof}
\begin{lemma}\label{lem: E(1)} The image of the $E(1)$--transfer is the
polynomial subalgebra generated by $h_{1,0}$.
\end{lemma}
\begin{proof}
From the previous lemma, we see that $f_{1}^{E(1)}(b_{m})$ is
nontrivial iff $m = 2^{t}-2$ for some $t \geq
1$. Note that if $t \geq 2$, then $m$ is even. So a monomial
consisting of only $b_{m}$ of the form above with at least one
index $m > 0$ cannot be a summand of a $P_{1}^{0}$--annihilated (hence
$E(1)$--annihilated) element.
\end{proof}
For $m >1$, even in the case of $E(2)$, we have not yet been able to
calculate the whole image of the $E(m)$--transfer, however, we have the
following result.
\begin{proposition}\label{prop: stem 1} For each $m \geq 1$, the subalgebra
$\f2 [h_{m,s}|0 \leq s <m]$ of $H^{*,*}(E(m)) = \f2 [h_{t,s}|s <m \leq
t] $ is in the image the $E(m)$--transfer.
\end{proposition}
\begin{proof}
Since the transfer is an algebra homomorphism, it suffices to show
that for each $0 \leq s \leq (m-1)$,
\begin{itemize}
\item $b_{2^{s} (2^{m}-1)-1} \in P_{E(m)}
      H_{*} (\z2)$, and
\item $f_{1}^{E(m)}$ sends $b_{2^{s} (2^{t}-1)-1}$ to
      $[\xi_{t}^{2^{s}}]$.
\end{itemize}
The first item is immediate from \fullref{lem: E(m)-annihilated}.
In fact, $b_{2^s (2^m -1) -1}$ is $E(n)$--annihilated for all $n >
m$. For the second item, it suffices to note that in a presentation
of $2^{s} (2^{t} -1)-1$ in the form \eqref{eq: nontrivial image},
\[
2^{s} (2^{t} -1)-1= 2^{s_{1}} (2^{t_{1}} -1) + \dotsb  +
2^{s_{\ell}}(2^{t_{\ell}}-1) -1,
\]
it follows that $2^{s} = 2^{s_{1}} + \dots  + 2^{s_{\ell}}$ by
\fullref{lem: image of Tnm}. This equation clearly has only one
possible solution: $\ell =1$ and $s_{1}=s$. But then $t_{1} =t$, so we
are done.
\end{proof}

Taking limits over all $m \geq 1$, we obtain the following result.
\begin{corollary}\label{cor: image of PD(m)}
For each $m \geq 1$, any element of the form 
$$h_{i_{0},\dotsc i_{m-1}}= h_{m,0}^{i_{0}}\dotsb h_{m,m-1}^{i_{m-1}},$$ 
where $i_{m-1} >0$ in the algebra $R$ of \fullref{thm: P inverse limit}
is contained in the image of the total
$\mathcal{E}$--transfer $\varphi^{\mathcal{E}}$.
\end{corollary}
\begin{proof}
First of all, observe that the family of maximal sub-Hopf algebras
$E(m)$ forms a cofinal system in $\mathcal{E}$. Thus to construct an
element in the inverse limit
$$\smash{\textstyle\invlim{E}}\bigoplus_{n}[P_{E}H_{*}(BV_{n})]_{GL(n)},$$
it suffices to define it on the $E(m)$s.

Now for any element of the form in the Corollary, we can define its
preimage as the compatible sequence which is zero for $E \neq E(m)$,
and
$$b_{i_{0}, \ldots i_{m-1}} = b_{2^{0}(2^{m}-1)-1}^{i_{0}}\ldots
b_{2^{m-1}(2^{m}-1)-1}^{i_{m-1}}$$
when $E=E(m)$. The fact that this
sequence is compatible is because the restriction of $b_{i_{0}, \dotsc
i_{m-1}}$ to any intersection $E(m) \cap E(m')$ is trivial. Here we
have used the condition that $i_{m-1} \geq 1$.
\end{proof}
\begin{corollary}\label{cor: The D transfer} The subalgebra of
$H^{*,*} (D)$ generated by $h_{t,t-1}$, $t \geq 1$, is contained in
the image of the $D$--transfer.
\end{corollary}
\begin{proof}
It is clear that $b_{2^{t-1} (2^{t}-1)-1}$ is
$D$--annihilated. On the other hand, we have that $\smash{f_{*} (b_{2^{t-1}
(2^{t}-1)-1}) = (f_{*} b_{2^{t}-2})^{2^{t-1}}}$. Since
$\smash{\xi_{t}^{2^{s}}} =0$
in $D_{*}$ for all $s \geq  t$, it follows easily that
the image of $\smash{b_{2^{t-1} (2^{t}-1)-1}}$ in $D_{*}$ is
$\smash{\xi_{t}^{2^{t}-1}}$ which is a cycle representing $h_{t,t-1}$.
\end{proof}

\begin{remark}
The set of elements of the form $h_{i_{0},\dotsc
i_{m-1}}$ with $i_{m-1} >0$ is the complete set of
$\st$--invariant monomials in $R$ (see Palmieri~\cite[page~433]{Palmieri99}).
Palmieri's theorem then predicts that some power, say $\alpha_{m}$, of
it comes from $H^{*,*} (\st )$. If our \fullref{conj: my own} is true, then
there would be another exponent $\beta_{m}$ for which
$\smash{b_{i_{0},\ldots i_{m-1}}^{\beta_{m}}}$ comes from a $GL$--coinvariant
element of $P_{\st} H_{*} V$. The two exponents are
probably not equal in general. For example, if $m=1$, then $h_{m,m-1}$
corresponds
to the element called $h_{0}$ in $H^{1,1} (\st)$ and $\alpha_{1}
= \beta_{1} =1$. But in the case $m=2$, then $\alpha_{2}=4$ where
$h_{2,1}^{4}$ corresponds to the element called $g$ in $H^{4,24}
(\st )$. However, $g$ is not an element of the $\st$--transfer
(see Bruner, H\`a and Hung~\cite{BHH}; also see the next section for a
quick proof). We do not know
what $\alpha_{m}$ and $\beta_{m}$ are in general. But for $m=2$, we
conjecture that one can take $\beta_{2}=6$ (Equivalently, this means
that the element denoted by $r \in H^{6,36} (\st )$ is in the image of
$\varphi_{6}^{\st}$.)
\end{remark}

One may wonder whether \fullref{prop: stem 1} describes all
elements in the image of the $E(m)$--transfer. This is not the case, as
our next example show.
\begin{proposition}\label{prop: an element in the image of the transfer}
$h_{3,0} h_{2,1}$ is in the image of the $\mathcal{E}$--transfer.
\end{proposition}
\begin{proof}
We will construct an element in the inverse limit
$\smash{\invlim{E}}[P_{E}H_{11}(BV_{2})]_{GL(2)}$ whose image in $R$
is $h_{3,0} h_{2,1}$.
Let $b \in H_{11} (BV_{2})$ be the sum $b = b_{6} b_5 + b_3 b_8 + b_9 b_2
+ b_{10} b_{1} + b_7 b_4$. By direct inspection, $b$ is
$E(2)$--annihilated.  We claim that $b$ represents a nontrivial
element in the coinvariant ring $[P_{E(2)} H_{11} (BV_{2})]_{GL(2)}$ but represents
a trivial class when replacing $E(2)$ by any other  $E \subset E(2)$.

Indeed, one verifies that the only $E(2)$--annihilated elements
in $H_{11} (BV_{2})$ are $b_{0}b_{11}$, $b$ and their obvious
permutations. It follows that $b$ represent a nontrivial element in
$[P_{E(2)} H_{11} (BV_{2})]_{GL(2)}$.

On the other hand, if $\sigma \in GL(2)$ denote the matrix
$\bigl(\begin{smallmatrix} 1&1 \\1 &0 \end{smallmatrix}\bigr)$ and
$\tau$ the standard permutation, then we have
\begin{align*}
b_{11} b_{0} + \sigma (b_{11} b_{0}) =& b + \tau  b + b_{0} b_{11} &
b_{9} b_{2} + \sigma (b_{9} b_{2}) =& b + b_{11} b_{0}.
\end{align*}
Thus if $E$ is a sub-algebra of $E(2)$ such that $b_{9}b_{2}$ is
$E$--annihilated, then $b$ represents the trivial class in the
corresponding coinvariant ring. This condition is true for the
subalgebras $E(2) \cap E(m)$ for any other $m$. We can now define an
element in $\smash{\invlim{E}} [P_{E} H_{11} (BV_{2})]_{GL(2)}$ by taking its
value to be $[b]$ when $E=E(2)$ and zero for other $E= E(m)$.

It remains to verify that the image of this element under the
$\mathcal{E}$--transfer is $h_{3,0}h_{2,1}$. For this purpose, note
that under the map $f_*\co H_{*} (BZ/2)\to
 \wwbar{\st}_{*}$, we see that $f_{*}(b_{3}) = \smash{\xi_{1}^{4}}$, $f_{*}(b_{9})
 = \smash{\xi_{1}^{10}+\xi_{1}^{4}\xi_2^{2}}$, $f_{*}(b_{1}) =
 \smash{\xi_{1}^{2}}$, and $f_{*}(b_{4}) =
\smash{\xi_{1}^{5}+\xi_{1}^{2}\xi_{2}}$,
 all are mapped to zero in $\smash{\overline{E(2)}^{*}}$. It follows that the
restriction of $b$ to ${(\overline{E(2)}^{*})^{2}}$ is $b_{6}
b_{5}$ which represents $h_{3,0} h_{2,1}$.

We remark that this element does not come from $[P_{D}
H_{11} (BV_{2})]_{GL(2)}$. For the only nontrivial element in the
latter is $\sum_{i=1}^{10}b_{i} b_{11-i}$ and its image in
$[P_{ E(2)}H_{11}(BV_{2})]_{GL(2)}$ is $b + \tau b  \equiv 0$. Our
\fullref{conj: my own} then predicts that $[b]$ is
nilpotent.
\end{proof}
Another useful question is to find a necessary criteria so that an
element in $H^{*,*}(E(m))$ is not in the image of the
$E(m)$--transfer. For this, we have the following.

\begin{proposition}\label{prop: same s}
Any element in $H^{*,*} (E(m))$ which contains a monomial of
the form $h = h_{t_{1}, s} \dotsb h_{t_{k},s}$ where not all $t_{i} =
m$, is not in the image of the $E(m)$--transfer.
\end{proposition}
\begin{proof}
By way of contradiction, suppose that $h$ is a nontrivial summand
of an element which is in the image of the $E(m)$--transfer. This
implies that there is an $E(m)$--annihilated element which contains
the monomial
\[
b = b_{2^{s} (2^{t_{1}}-1)-1} \dotsb b_{2^{s} (2^{t_{k}}-1)-1},
\]
as a nontrivial summand. Equivalently, the dual of $b$
\[
x = x_{1}^{2^{s} (2^{t_{1}}-1)-1} \dotsb x_{k}^{2^{s} (2^{t_{k}}-1)-1},
\]
is $E(m)$--indecomposable. We will show that this last statement is not
true. Indeed, observe that for each $t$, $2^{s} \notin 2^{s}
(2^{t}-1)-1 = (2^{s+t} -1) - 2^{s}$. Thus by \fullref{lem: Pst action},
$x^{2^{s} (2^{t}-1) -1}$ is $P_{\ell}^{s}$--annihilated for
any $\ell$. It follows that if there exist $t_{i} > n$, say $t_{1}$,
then $P_{t_{1}-1}^{s} \in E(m)$ and
\[
P_{t_{1}-1}^{s} (x_{1}^{2^{s+t_{1} -1}} x_{2}^{2^{s}(2^{t_{2}}-1)-1}\dotsb
x_{k}^{2^{s}(2^{t_{k}}-1)-1}) =
x_{1}^{2^{s} (2^{t_{1}}-1)-1} \dotsb x_{k}^{2^{s} (2^{t_{k}}-1)-1}.
\]
The lemma is proved.
\end{proof}
\begin{example}\label{ex: Palmieri example}
Here is an example of the usefulness of the last proposition. In
\cite{Palmieri99}, Palmieri asserted that there is an $\st$--invariant
element of $R$ of the form
\[
z^{12,80} =  h_{2,0}^{8} h_{3,1}^{4} + h_{3,0}^{8} h_{2,1}^{4} +
h_{2,1}^{11} h_{3,1}.
\]
This element contains the monomial $\smash{h_{2,1}^{11} h_{3,1}}$ of the form
described in \fullref{prop: same s}, thus $z^{12,80}$ or any
of its powers cannot be in the image of the $E(2)$--transfer.
\end{example}

\section{Applications to the rank $4$ transfer}\label{sec:info}
The cohomology of the Steenrod algebra in cohomological degree $4$
has been completely determined by Lin and Mahowald \cite{Lin-Mahowald}. In particular,
there are three generators, namely $d_{0} \in H^{4,18}(\st )$,$e_{0}
\in H^{4,21}(\st )$, and $g \in H^{4,24}(\st )$ whose restriction
to the sub-Hopf algebra $E(2)$ are nontrivial. According
to Zachariou \cite{Zachariou67,Zachariou73}, their images in
$H^{4,*}(E(2))$ are $\smash{h_{2,0}^{2}h_{2,1}^{2}}$,
$\smash{h_{2,0}h_{2,1}^{3}}$ and $\smash{h_{2,1}^{4}}$
respectively. Furthermore, $\smash{i^{\st}_{E(2)}}$ is a monomorphism
in these degrees.
In this section, we prove \fullref{thm:1.1}. That $g$, in fact the whole family of generators $g_i =
(Sq^{0})^{i} g$, are not in the image of the rank $4$ transfer
was first proved by Bruner, H\`a and Hung~\cite{BHH}. We provide here
a different proof
which is less calculational. The elements $d_0$ and $e_0$ are in Hung's
list of conjectural elements that must be in the image of the
transfer, provided the Singer conjecture holds. Our results thus
partially complete his list\footnote{Recently, Nam
\cite{Nam-tr4}  claimed to have verified most cases in Hung's
list by using a completely different method.}.
Recall first that under the map $\smash{f_{4}^{E(2)}}$, the images of
$b_{k}$ for $k \leq 20$ are
all trivial, except the following:
{\small
$$\begin{array}{lclclclcl}
b_{2}\longmapsto\xi_{2} & \qua & b_{5}\longmapsto\xi_{2}^{2} & \qua &
  b_{6}\longmapsto\xi_{3} & \qua & b_{8}\longmapsto\xi_{2}^{3} & \qua & 
  b_{11}\longmapsto\xi_{2}^{4} \\
b_{12}\longmapsto\xi_{2}^{2}\xi_{3} & \qua & b_{13}\longmapsto\xi_{3}^{2} &
  \qua & b_{14}\longmapsto\xi_{4} & \qua & 
  b_{16}\longmapsto\xi_{2}\xi_{3}^{2} & \qua &
  b_{20}\longmapsto\xi_{2}\xi_{4}+\xi_{3}^{3}
\end{array}$$}%
The proof of \fullref{thm:1.1} is divided
into three parts, corresponding to the three generators in question.
\subsection{$g$ is not in the image of the transfer}
We prove by contradiction. Suppose that there is an element $z \in
P_{\st} H_{20} (BV_{4})$ such that $\varphi_{4}^{\st} ([z]) = g$. Since the restriction of $g$ to $H^{*}
(E(2))$ is $h_{2,1}^{4}$, it follows  easily that $z$ must contains the monomial $b_{5}
b_{5}b_{5}b_{5}$ as a nontrivial summand. Dually, we have that the monomial
$x_{1}^{5}x_{2}^{5}x_{3}^{5}x_{4}^{5} =5555$ is indecomposable in
the $\st$--module $H^{20} (BV_{4})$.  We will show that this latter statement is not true.
In \cite{Wood-Peterson}, Wood proved Peterson's conjecture which says
that a monomial in degree $d$ with exactly $r$ odd exponents such that
$\alpha (d+r) >r$ is $\st$--decomposable; where $\alpha (m)$
counts the number of non-zero digits in the binary expansion of $m$.
His proof makes use of an important observation, nowadays known as Wood's
$\chi$--trick, that for any two homogeneous polynomials
$u$ and $v$ in the polynomial algebra $H^{*} (BV_{n})$ and any
Steenrod operation $\theta$, $u (\theta v)$ is $\st$--decomposable iff $(\chi (\theta) u) v$
is, where $\chi$ is the canonical conjugation in $\st$.

Writing $u \equiv v$ whenever $u-v$ is $\st$--decomposable, we have
\begin{align*}
5555 = Sq^{8} (2222) \times 1111 \equiv & 2222 \times (\chi Sq^{8}) 1111\\
\equiv & 2222 \times [(4422) + (8211)]\\
\equiv & (6,6,4,4) + (10,4,3,3) \equiv (10,4,3,3),
\end{align*}
where a monomial in brackets means that we take the sum of all possible
permutations of that monomial. But the monomial $(10,4,3,3)$ is
$\st$--decomposable because $\alpha (20 + 2) = 3 > 2$. We have a
contradiction.

It can be easily  seen that the proof above also works for any other
generator $g_{i} = (\smash{\widetilde{Sq}^{0}})^{i} g$.

\subsection{$d_{0}$ is in the image of the transfer}\label{subsec: d0}
We will show that there exists an element $z \in P_{\st} H_{14} BV_{4}$
such that it contains an odd number of permutations of $2255$.
Assuming that such an element exists, then its image under the
canonical maps
\[
P_{\st} H_{14} (BV_{4}) \hookrightarrow  P_{E(2)} H_{14} (BV_{4}) \to
[P_{E(2)} H_{14} (BV_{4})]_{GL(4)},
\]
is the equivalent class of $2525$ because in degree $14$, there are
only two possible type of monomials, namely $2255$ and $2228$ that
maps nontrivially to $\overline{E(2)}^{*}$. The latter monomial obviously cannot be
a nontrivial summand of any $\st$--annihilated element. It follows
that the image of $z$ under the $\st$--transfer is an element in
$H^{4,14} (\st )$ whose restriction to $H^{4,14} (E(2))$ is
$h_{2,0}^{2} h_{2,1}^{2}$. Thus $z$ is a chain-level representation of $d_{0}$.

In fact, one can verify that the $\st$--annihilated element
$$
z = x + (2,3) x + (1,3) x + (3155+5513+5135+5315+5333),
$$
where 
$$
x = (2255+2165+1256+1166+4253+4163+3263+2435+1436+2336+4433), 
$$
is a representation of $d_{0}$.
Indeed, to verify that this sum is indeed $\st$--annihilated, we need only consider the effects of $Sq^{1}, Sq^{2}$ and $Sq^{4}$ because of the unstable condition. By direct calculation, we have
\begin{align*}
Sq^1 x = Sq^1 4433 =& 4333+3433,\\
Sq^2 x =&  3153 + 1335 + 3333,\\
Sq^4 x =&  1333 + 3133. 
\end{align*}
It follows that $x + (2,3)x + (1,3)x$ is $Sq^1$-- and $Sq^4$--nil. Moreover, 
$$
Sq^2 [x + (2,3)x + (1,3)x ] = 3153 + 3513 + 3315 + 5133 + 3333.
$$
The extra summand $(3155+5513+5135+5315+5333)$ is needed to kill off the effect of $Sq^2$.

\subsection{$e_{0}$ is in the image of the transfer}
There is only one type of monomial, namely $2555$, in $H_{17} (BV_{4})$
whose image in
$(\overline{E(2)}^{*})^{4}$ is nontrivial. As in the proof for
$d_{0}$, it suffices to show that there exists an element of $P_{\st}
H_{17} BV_{4}$ which contains an odd number of permutations of $2555$.

We exhibit an explicit element of such form below.
\begin{align*}
2555&+1655+18(53)+17(63)+14(75)+13(76)+14(93)+23(93)\\
&+12(95)+11,10,5+1169+12(11,3)+4(355)+11,12,3+114,11+\\
&+1187+2177+112,13+111,14+3356+3635+3563\\
&+5336+5633+5363 +6(335)+8333+7433+7253+7163\\
&+2933+1,10,33+2735+2375+2357+1736+1376+1367.
\end{align*}



\bibliographystyle{gtart}
\bibliography{link}

\end{document}